\documentclass[a4paper,12pt]{amsart}

\usepackage{amsthm, amsmath, latexsym, amsfonts, epsfig}
\usepackage[all]{xy}
\usepackage{amssymb}

\newtheorem{Lemma}{Lemma}[section]
\newtheorem{Proposition}[Lemma]{Proposition}
\newtheorem{Theorem}[Lemma]{Theorem}

\newtheorem{Corollary}[Lemma]{Corollary}

\newcommand{\Pic}{\mathrm{Pic}}

\newcommand{\irr}{\mathrm{irr}}
\newcommand{\Z}{\mathbb{Z}}

\newcommand{\Q}{\mathbb{Q}}

\newcommand{\bP}{\mathbb{P}}

\newcommand{\sgbar}{\overline{S}_g}

\newcommand{\sgnm}{\overline{S}_{g,n}^{\hspace{0.05cm}(m_1, \ldots, m_n)}}

\title[Cohomology of moduli of curves with level structures]
{On the rational cohomology of \\ moduli spaces of curves \\ with
level structures}
\author{Gilberto Bini and Claudio Fontanari}

\email{gilberto.bini@unimi.it} \curraddr{{\sc Dipartimento di
Matematica \\ Universit\`a degli Studi di Milano \\ Via C. Saldini 50 \\
20133 Milano \\ Italy.}}

\email{fontanar@science.unitn.it}\curraddr{
{\sc Dipartimento di Matematica \\  Universit\`a degli Studi di Trento\\
Via Sommarive 14 \\ 38123 Trento \\ Italy.}}


\begin{document}

\begin{abstract}
We investigate low degree rational cohomology groups of smooth
compactifications of moduli spaces of curves with level structures.
In particular, we determine $H^k(\sgbar, \Q)$ for $g \ge 2$ and $k
\le 3$, where $\sgbar$ denotes the moduli space of spin curves of
genus~$g$.
\end{abstract}

\maketitle

\section{Introduction}\label{sec1}

Moduli spaces $M^\Gamma$ of curves with level structures are
obtained by taking the quotient of the Teichm\"uller space by a
finite index subgroup $\Gamma$ of the mapping class group. As such,
they yield natural coverings of the moduli space $M_{g,n}$ of smooth
genus $g$ curves with $n$ marked points.

The geometry of $M^\Gamma$ turns out to be better understood when
$\Gamma$ contains the Torelli group, as explained in \cite{Hain}.
For instance, $\Gamma$ can be a spin mapping class group, so that
$M^\Gamma$ is the moduli space $S_g^\pm$ parameterizing pairs
(smooth genus $g$ curve $C$, even/odd theta characteristic on $C$).
Under this assumption, a recent result by Putnam \cite{Put} shows
that $H^2(M^\Gamma, \Q) \cong \Q$ for $g \ge 5$.

Here instead we investigate the rational cohomology of the canonical
compactification $\overline{M}^\Gamma$ of $M^\Gamma$ over the moduli
space $\overline{M}_{g,n}$ of stable curves. Indeed, the description
of the boundary of $\overline{M}^\Gamma$ provided by \cite{Bog} and
\cite{Bogg} allows us to adapt the inductive approach introduced for
$\overline{M}_{g,n}$ in \cite{ArbCor:98} and recently refined in
\cite{ArbCor:08}. In particular, we are able to show that
$H^1(\overline{M}^\Gamma, \Q)=0$ (see Theorem~\ref{acca1}) and to
determine a set of free generators for $\Pic(\overline{M}^\Gamma)
\otimes \Q$ (see Theorem~\ref{pic}).

In the special case of spin moduli spaces, where a geometrically
meaningful compactification has been constructed in \cite{Cor}, we
obtain even stronger results. Namely, we get the vanishing of the
third cohomology group (see Theorem~\ref{spinvanishing}) and a
complete description of the second cohomology group (see
Corollary~\ref{h2}).

Throughout,  we work over the field ${\mathbb C}$ of complex numbers
and all cohomology groups are intended to be with rational
coefficients.

The first named author has been partially supported by "FIRST" Universit\`{a} di Milano and by MIUR Cofin 2008 - Variet\`{a} algebriche: geometria, aritmetica, strutture di Hodge. The second named author has been partially supported by MIUR Cofin 2008 - Geometria delle variet\`{a} algebriche e dei loro spazi di moduli.

\section{The first Betti number of moduli spaces of curves with level structures}\label{sec2}

First, we recall some notation and basic facts on level structures:
for more details, the reader is referred, for instance, to
\cite{Bog}.

Let $\Sigma_{g,n}$ be a compact genus $g$ surface with $n$ marked
points and let $\Gamma_{g,n}$ be its mapping class group. Denote by
$T_{g,n}$ the Teichm\"{u}ller space and by $M_{g,n}$ the moduli
space $T_{g,n}/\Gamma_{g,n}$. A {\em level} $\Gamma$ is a subgroup
of $\Gamma_{g,n}$. It is finite if $\Gamma$ has finite index in
$\Gamma_{g,n}$. It is {\em Galois} if $\Gamma$ is normal in
$\Gamma_{g,n}$. The functor of curves with $\Gamma$-level is
represented by the analytic stack $[T_{g,n}/\Gamma]$, which is
called a level structure over the moduli space $M_{g,n}$. Clearly, a
level structure is a finite connected covering of $M_{g,n}$.

In what follows, we will focus on some particular level structures.
Fix $g \geq 2$ and a basis of the first homology group of
$\Sigma_{g,n}$ so that the intersection form is given by the $2g
\times 2g$ matrix
$$
\left(
\begin{array}{cc}
0 & I \\
-I & 0
\end{array}
\right).
$$

There exists a surjective homomorphism $\Gamma_{g,n} \rightarrow
Sp(2g, {\mathbb Z})$, where $Sp(2g, {\mathbb Z})$ is the symplectic
group of $2g \times 2g$ matrices with integer entries. The kernel of
this homomorphism is called the Torelli group ${\mathcal T}_{g,n}$.
In some applications, we will take into account finite levels that contain
the Torelli group. We briefly review some examples which fit into
this picture.

{\em Example A.} For any integer $m \geq 2$, consider the surjective
homomorphism $\Gamma_{g,n} \rightarrow Sp(2g, {\mathbb Z}/{m \mathbb
Z})$ which maps an element $\gamma \in \Gamma_{g,n}$ to the
homomorphism induced by $\gamma$ on the homology of $\Sigma_{g,n}$
mod $m$. The kernel of this homomorphism is called {\em the
Abelian level of order m}. For short, we will denote this level by
$(m)$. By definition, it contains the Torelli group.

Those of Example A are particular examples of geometric levels. Indeed, if
$\Pi^{\lambda}$ is a subgroup of the fundamental group $\Pi_{g,n}$ of $\Sigma_{g,n}$,
then the {\em geometric level} determined by $\Pi^{\lambda}$ is the kernel
$\Gamma^{\lambda}:=ker\rho_{\lambda}$ of the natural representation:
$$ \rho_{\lambda}: \Gamma_{g,n} \rightarrow Out(\Pi_{g,n}/\Pi^{\lambda})).
$$

A level $\Gamma^\lambda$ is said to be {\em fine} if $\Gamma^\lambda
\subset (m)$. In geometric terms, this means that ${M}^{\lambda}$ is
a connected \'{e}tale covering of $M^{(m)}$.

{\em Example B.} Let $V= H_1(\Sigma_{g,n}, {\mathbb Z}/2{\mathbb
Z})$ and let $\cdot$ denote the ${\mathbb Z}/{2 \mathbb Z}$-
intersection form on $V$. A ${\mathbb Z}/{2 \mathbb Z}$-quadratic
form on $\Sigma_{g,n}$ is a function $Q: V \rightarrow {\mathbb
Z}/{2\mathbb Z}$ such that
$$
Q(x+y)=Q(x)+Q(y)+x\cdot y
$$
for all $x,y \in V$. The isomorphism class of $Q$ is determined by
its Arf invariant. Given $Q$, let $G(Q)$ denote the subgroup of
$\Gamma_{g,n}$ which preserves $Q$. Since there are two isomorphism
classes of quadratic forms, there are two such groups, which are
called the spin mapping class groups. It is well known that they
contain the Torelli group. The space $T_{g,n}/G(Q)$ is
the union of two connected components, $S_{g,n}^{+} \sqcup
S_{g,n}^{-}$, which are called the moduli spaces of even (respectively,
odd) spin structures.

A canonical compactification $\overline{M}^{\lambda}$ of
$M^{\lambda}$ is obtained by taking the normalization of the
Deligne-Mumford compactification $\overline{M}_{g,n}$ in the
function field of $M^{\lambda}$. As proved in Proposition 1.6.8 in
\cite{Bogg}, if $\Gamma^{\lambda}$ is a finite fine level, then
${\overline M}^{\lambda}$ is represented by a projective variety.

Finally, we recall the geometry of the boundary as it is described
in \cite{Bog}, Proposition 2.2 and Theorem 2.3, and in \cite{Bogg},
Theorem 2.5.1 and Theorem 2.7.4.
Let $\lambda$ be a fine geometric level structure over $M_{g,n}$
such that its compactification ${\overline M}^{\lambda}$ is smooth.
An irreducible component of its Deligne-Mumford boundary,
corresponding to reducible curves, is isomorphic to
$$
{\overline M}^{\lambda_1}_{g_1, n_1+1} \times {\overline M}^{\lambda_2}_{g_2, n_2+1}
$$
for some $g_1+g_2=g$, $n_1+n_2=n$ and with $\lambda_i$ suitably
defined geometric levels (see \cite{Bogg}, p. 25). Analogously, the
closure of each stratum parameterizing singular irreducible curves
is isomorphic to
$$
{\overline M}_{g-h,n+2h}^{\lambda_{\sigma}},
$$
where $h \ge 1$ and $\lambda_{\sigma}$ is a suitably defined
geometric level (see \cite{Bogg}, p. 35).

In the special case of the spin mapping class group, a modular
compactification \`{a} la Deligne-Mumford has been constructed by
Cornalba in \cite{Cor}, where one can find an explicit description
of the boundary as well (see \cite{Cor}, \S~7).

Next, we turn to cohomological computations. Let $\Gamma$ be a level
of $\Gamma_{g,n}$. We point out that, analogously to $M_{g,n}$, the
homology of $M^{\Gamma}$ vanishes in high degree. In fact, the
following holds.

\begin{Proposition}\label{harer} Let $\Gamma \subset \Gamma_{g,n}$. Then
$H_k(M^{\Gamma}, {\mathbb Q})=0$ for $k > c(g,n)$, where
$$
c(g,n)= \left\{
\begin{array}{cc}
n-3 & g=0; \\
4g-5 & g > 0, n=0; \\
4g-4+n & g>0 , n > 0. \end{array}
\right.
$$
\end{Proposition}

\proof As proved in \cite{Harer:86}, Theorem 1.3, for every $n \ge 1$ there exists a
$\Gamma_{g,n}$-equivariant homeomorphism of $T_{g,n}$ onto the arc complex. A fortiori,
this is a $\Gamma$-equivariant homeomorphism. Thus, the standard proof of Harer's
vanishing theorem \cite{Harer:86} for the high degree homology of $M_{g,n}$
(see for instance \cite{ArbCor:08}, Lemma 2) shows that the same result holds
for any such $M^{\Gamma}$ as well. On the other hand, the case $n = 0$ is ruled out
by a standard spectral sequence argument (we refer again to \cite{ArbCor:08},
last paragraph of Section 5, which adapts verbatim to our context).

\qed

As a consequence, we show how  the inductive approach of
\cite{ArbCor:98} to the rational cohomology of $\overline{M}_{g,n}$
applies verbatim also to level structures and reduces the
computation for a fixed degree to a few initial cases. In
particular, for degree $1$ the inductive basis is provided by the
following easy fact.

\begin{Lemma}\label{basis}
Denote by ${\overline M}^{\lambda}$ the smooth canonical compactification
of a level structure of $\overline{M}_{0,4}$ or $\overline{M}_{1,1}$.
Then ${\overline
 M}^{\lambda}$ is isomorphic to $\mathbb{P}^1$.
\end{Lemma}

\proof The degree $d$ covering map can be ramified only over points
parameterizing curves with non-trivial automorphisms, whose number
is $0$ for $\overline{M}_{0,4}$ and $2$ for $\overline{M}_{1,1}$
(see, for instance, \cite{Har}, Corollary IV.4.7). Hence Hurwitz
formula for the genus $g$ of ${\overline M}^{\Gamma}$ yields $2g-2 =
d(-2)$ and $2g-2 \le d(-2)+2(d-1)$, respectively; in both cases it
follows that $g = 0$.

\qed

>From Lemma \ref{basis}, Proposition \ref{harer}, Poincar\'{e}
duality and the long exact sequence in cohomology with compact
support
$$
\ldots \to H^k_c(M^\lambda) \to H^k({\overline M}^{\lambda}) \to
H^k(\partial {\overline M}^{\lambda}) \to \ldots
$$
we obtain the following result (see for instance \cite{ArbCor:08},
proof of Corollary 1 to Theorem 6).

\begin{Theorem}
\label{acca1} Let $\lambda$ be a fine geometric level structure over $M_{g,n}$
such that the canonical compactification ${\overline M}^{\lambda}$ over
$\overline{M}_{g,n}$ is smooth. Then $H^1({\overline M}^{\lambda}, \mathbb{Q})=0$.
\end{Theorem}

\section{Low degree cohomology of spin moduli spaces}

In the previous sections we have already mentioned the moduli space
of spin curves constructed in \cite{Cor}. More generally, for all integers
$g$, $n$, $m_1, \ldots, m_n$, such that $2g-2+n>0$, $0 \le m_i \le 1$
for every $i$, and $\sum_{i=1}^n m_i$ is even, one can consider the moduli
spaces
\begin{eqnarray*}
\sgnm &:=& \{ [(C, p_1, \ldots, p_n; \zeta; \alpha)]:
(C, p_1, \ldots, p_n) \textrm{ is a genus $g$} \\
& & \textrm{quasi-stable projective curve with $n$ marked points}; \\
& &\zeta \textrm{ is a line bundle of degree $g-1+ \frac{1}{2} \sum_{i=1}^n m_i$
on $C$} \\
& & \textrm{having degree $1$ on every exceptional component of $C$,} \\
& &\textrm{and } \alpha: \zeta^{\otimes 2} \to \omega_C(\sum_{i=1}^n
m_i p_i)
\textrm{ is a homomorphism which} \\
& & \textrm{is not zero at a general point of every
non-exceptional}\\
& & \textrm{component of $C$} \}.
\end{eqnarray*}

Here we prove the following result on the rational cohomology of $\sgnm$:

\begin{Theorem}\label{spinvanishing}
For every $g$, $n$ and $(m_1, \ldots, m_n)$ as above, we have
$$
H^1(\sgnm, \Q) = H^3(\sgnm, \Q) = 0.
$$
\end{Theorem}

We are going to apply the inductive strategy developed by Arbarello
and Cornalba in \cite{ArbCor:98} for the moduli space of curves.
Namely, we consider the long exact sequence of cohomology with
compact supports:
\begin{equation}\label{exact}
\ldots \to H^k_c(S_{g,n}^{\hspace{0.05cm}(m_1, \ldots, m_n)}) \to
H^k(\sgnm) \to H^k(\partial S_{g,n}^{\hspace{0.05cm}(m_1, \ldots, m_n)})
\to \ldots
\end{equation}
Hence, whenever $H^k_c(S_{g,n}^{\hspace{0.05cm}(m_1, \ldots, m_n)})
= 0$, there is an injection $H^k(\sgnm)$ $\hookrightarrow
H^k(\partial S_{g,n}^{\hspace{0.05cm}(m_1, \ldots, m_n)})$.
Moreover, from \cite{Cor}, \S~3, it follows that each irreducible
component of the boundary of $\sgnm$ is the image of a morphism:
$$
\mu_i: X_i \to \sgnm
$$
where either
$$
X_i = \overline{S}_{a, s+1}^{\hspace{0.05cm} (u_1, \ldots, u_{s+1})}
\times \overline{S}_{b, t+1}^{\hspace{0.05cm} (v_1, \ldots, v_{t+1})}
$$
with $a+b=g$, $s+t=n$, and $\sum_{i=1}^s u_i +\sum_{i=1}^t v_i =
\sum_{i=1}^n m_i$; or
$$
X_i = \overline{S}_{g-1, n+2}^{\hspace{0.05cm} (m_1, \ldots, m_n,
m_{n+1}, m_{n+2})}.
$$
Finally, exactly as in \cite{ArbCor:98}, Lemma~2.6, a bit of Hodge
theory implies that the map
\begin{equation}
H^k(\sgnm) \to \oplus_i H^k(X_i) \label{injective}
\end{equation}
is injective whenever $H^k(\sgnm) \to H^k(\partial
S_{g,n}^{\hspace{0.05cm}(m_1, \ldots, m_n)})$ is. So we obtain the
claim of Theorem~\ref{spinvanishing} by induction, provided we show
that $H^1_c(S_{g,n}^{\hspace{0.05cm}(m_1, \ldots, m_n)}) =
H^3_c(S_{g,n}^{\hspace{0.05cm}(m_1, \ldots, m_n)}) = 0$ for almost
all the values of $g$, $n$, and we check that $H^1(\sgnm) =
H^3(\sgnm) = 0$ for all the remaining values of $g$, $n$.

>From Proposition~\ref{harer} and Poincar\'e duality we deduce that
$$
H^1_c (S_{g,n}^{\hspace{0.05cm}(m_1, \ldots, m_n)}) = 0
$$
for any $g \ge 2$, for $g=1$, $n \ge 2$, and for $g=0$, $n \ge 5$;
and
$$
H^3_c (S_{g,n}^{\hspace{0.05cm}(m_1, \ldots, m_n)}) = 0
$$
for any $g \ge 3$, for $g=2$, $n \ge 2$, for $g=1$, $n \ge 4$, and
for $g=0$, $n \ge 7$.
Hence we have only to verify that
\begin{eqnarray}
\label{eq_one} H^1(\overline{S}_{0,n}^{\hspace{0.05cm}(m_1, \ldots, m_n)})
&=& 0,
\hspace{0.2cm} n \le 4 \\
\label{eq_two} H^3(\overline{S}_{0,n}^{\hspace{0.05cm}(m_1, \ldots, m_n)})
&=& 0, \hspace{0.2cm} n \le 6 \\
\label{eq_three} H^1(\overline{S}_{1,1}^{\hspace{0.05cm}(m)}) &=& 0 \\
\label{eq_four} H^3(\overline{S}_{2}) &=& 0 \\
\label{eq_five} H^3(\overline{S}_{2,1}^{\hspace{0.05cm}(m)}) &=& 0 \\
\label{eq_six}  H^3(\overline{S}_{1,n}^{\hspace{0.05cm}(m_1, \ldots, m_n)})
&=& 0, \hspace{0.2cm} n \le 3.
\end{eqnarray}

The first two checks are straightforward: since any two divisors of the same
degree on $\bP^1$ are linearly equivalent, there are natural isomorphisms
\begin{equation}
\label{zero}
 \overline{S}_{0,n}^{\hspace{0.05cm}(m_1, \ldots, m_n)}
\cong \overline{M}_{0,n}.
\end{equation}

Since $H^k(\overline{M}_{0,n})=0$ for every odd $k$ by Keel's
results (see \cite{Keel:92}), (\ref{eq_one}) and (\ref{eq_two})
easily follow.

As for (\ref{eq_three}), first of all notice that
$\overline{S}_{1,1}^{\hspace{0.05cm}(1)} = \emptyset$
by degree reasons.
Next, recall that $\overline{S}_{1,1}^{\hspace{0.05cm}(0)}$ is
the union of two connected components
$\overline{S}_{1,1}^{\hspace{0.05cm}(0), +}$ and
$\overline{S}_{1,1}^{\hspace{0.05cm}(0), -}$,
corresponding respectively to even and odd spin structures.
If $(E; q_1)$ is a smooth $1$-pointed elliptic curve, the linear series
$\vert 2 q_1 \vert$ realizes $E$ as a two-sheeted covering of $\bP^1$
branched at $q_1$ and at other three points $q_2$, $q_3$, and $q_4$.
The curve $E$ carries one odd theta-characteristic ($L = \mathcal{O}_E$),
and three even ones (namely, $\mathcal{O}_E(q_1-q_2)$,
$\mathcal{O}_E(q_1-q_3)$, and $\mathcal{O}_E(q_1-q_4)$).
The uniqueness of the odd theta-characteristic on $E$ implies the
existence of a natural isomorphism
$$
\overline{S}_{1,1}^{\hspace{0.05cm}(0), -} \cong \overline{M}_{1,1}
\cong \bP^1,
$$
hence $H^1(\overline{S}_{1,1}^{\hspace{0.05cm}(0), -}) = H^1(\bP^1)=0$.
Finally we turn to $\overline{S}_{1,1}^{\hspace{0.05cm}(0), +}$.
We claim that there is a surjective morphism
$$
f: \overline{M}_{0,4} \longrightarrow
\overline{S}_{1,1}^{\hspace{0.05cm}(0), +}.
$$
Indeed, let $(C; p_1, p_2, p_3, p_4)$ be a $4$-pointed stable genus zero
curve. The morphism $f$ associates to it the admissible covering $E$
of $C$ branched at the $p_i$'s, pointed at $q_1$ and equipped with the
line bundle $\mathcal{O}_E(q_1 - q_2)$, where $q_i$ denotes the point of
$E$ lying above $p_i$. It follows that
$$
H^1(\overline{S}_{1,1}^{\hspace{0.05cm}(0), +}) \hookrightarrow
H^1(\overline{M}_{0,4}) = H^1(\bP^1) = 0
$$
and (\ref{eq_three}) is completely proved.

The proofs of (\ref{eq_four}) and (\ref{eq_five}) are similar.
Again, $\overline{S}_{2,1}^{\hspace{0.05cm}(1)} = \emptyset$ by degree
reasons and $\overline{S}_{2,n}^{\hspace{0.05cm}(0, \ldots,0)}$ is
the disjoint union of $\overline{S}_{2,n}^{\hspace{0.05cm}(0, \ldots, 0), +}$
and $\overline{S}_{2,n}^{\hspace{0.05cm}(0, \ldots, 0), -}$
Moreover, if $C$ is a smooth hyperelliptic curve and $q_i$ ($i=1, \ldots, 6$)
are the ramification points of the hyperelliptic involution, then $C$ carries
six odd theta-characteristics (namely, $\mathcal{O}_C(q_i)$, $i=1, \ldots,
6$) and ten even ones (namely, $\mathcal{O}_C(q_i+q_j-q_k)$, with $i$, $j$
and $k$ distinct). We claim that there are surjective morphisms:
\begin{eqnarray*}
f^+: \overline{M}_{0,6} &\longrightarrow& \overline{S}_{2}^+ \\
f^-: \overline{M}_{0,6} &\longrightarrow&  \overline{S}_{2}^- \\
g^+: \overline{M}_{0,7} &\longrightarrow&
\overline{S}_{2,1}^{\hspace{0.05cm}(0,0), +} \\
g^-: \overline{M}_{0,7} &\longrightarrow&
\overline{S}_{2,1}^{\hspace{0.05cm}(0,0), -}
\end{eqnarray*}
In order to define $f^+$ and $f^-$, let $(C; p_1, \ldots p_6)$ be a
$6$-pointed, stable, genus zero curve. The morphism $f^+$ (respectively,
$f^-$) associates to it the admissible covering $Y$ of $C$ branched at the
$p_i$'s and equipped with the line bundle $\mathcal{O}_Y(q_1+q_2-q_3)$
(respectively, $\mathcal{O}_Y(q_1)$),
where $q_i$ denotes the point of $E$ lying above $p_i$.
As for $g^+$ and $g^-$, let $(C; p_1, \ldots p_7)$ be a
$7$-pointed stable genus zero curve. The morphism $g^+$ (respectively,
$g^-$) associates to it the admissible covering $Y$ of $C$ branched at the
$p_i$'s, pointed at one of the two points lying above $p_7$ (of course
different choices produce isomorphic curves) and equipped with the line bundle
$\mathcal{O}_Y(q_1+q_2-q_3)$ (respectively, $\mathcal{O}_Y(q_1)$),
where $q_i$ denotes the point of $E$ lying above $p_i$.
Hence we obtain injective maps in cohomology
$H^k(\overline{S}_{2}^+)\hookrightarrow H^k(\overline{M}_{0,6})$,
$H^k(\overline{S}_{2}^-)\hookrightarrow H^k(\overline{M}_{0,6})$,
$H^k(\overline{S}_{2,1}^{\hspace{0.05cm}(0,0), +}) \hookrightarrow H^k(\overline{M}_{0,7})$,
and $H^k(\overline{S}_{2,1}^{\hspace{0.05cm}(0,0), -}) \hookrightarrow H^k(\overline{M}_{0,7})$,
which reduce (\ref{eq_four}) and (\ref{eq_five}) to Keel's results mentioned above.

The proof of (\ref{eq_six}) turns out to be more involved. The case
$m_1=m_2=m_3=0$ has already been addressed in \cite{BF:09},
Lemma 4, so here we directly turn to $\overline{S}_{1,3}^{\hspace{0.05cm}(1,1,0)}$.
However, our inductive approach requires to handle $\overline{S}_{1,2}^{\hspace{0.05cm}(1,1)}$ too.

The boundary components of
$\overline{S}_{1,2}^{\hspace{0.05cm}(1,1)}$ are the following:
\begin{itemize}
\item $A_\irr$, whose general member is obtained from a smooth
$4$-pointed rational curve $C$ carrying the line bundle
$\mathcal{O}_C(P)$ by collapsing two marked points in
an ordinary node;
\item $B_\irr$, whose general member is obtained from a smooth
$4$-pointed rational curve carrying the line bundle $\mathcal{O}_C$
by joining two marked points with an exceptional component;
\item $A_{1, \emptyset}$, whose general member is obtained
by joining with an exceptional component
a smooth $1$-pointed elliptic curve $(E, p)$ carrying an even root of
$\mathcal{O}_E$ and a smooth $3$-pointed
rational curve $(C, q, p_1, p_2)$ carrying the line bundle $\mathcal{O}_C$;
\item $B_{1, \emptyset}$, whose general member is obtained by
joining with an exceptional component a smooth $1$-pointed elliptic
curve $(E,p)$ carrying the line bundle $\mathcal{O}_E$ and a smooth $3$-pointed
rational curve $(C, q, p_1, p_2)$ carrying the line bundle $\mathcal{O}_C$.
\end{itemize}

Next, we list the boundary components of
$\overline{S}_{1,3}^{\hspace{0.05cm}(1,1,0)}$:
\begin{itemize}
\item $A_\irr$, whose general member is obtained from a smooth
$5$-pointed rational curve $C$ carrying the line bundle
$\mathcal{O}_C(P)$ by collapsing two marked points in
an ordinary node;
\item $B_\irr$, whose general member is obtained from a smooth
$5$-pointed rational curve carrying the line bundle $\mathcal{O}_C$
by joining two marked points with an exceptional component;
\item $A_{1, \emptyset}$, whose general member is obtained by
joining with an exceptional component a smooth $1$-pointed elliptic curve $(E,p)$
carrying an even root of $\mathcal{O}_E$ and a smooth $4$-pointed rational
curve $(C, q, p_1, p_2, p_3)$ carrying the line bundle $\mathcal{O}_C$;
\item $B_{1, \emptyset}$, whose general member is obtained by
joining with an exceptional component a smooth $1$-pointed elliptic curve $(E,p)$
carrying the line bundle $\mathcal{O}_E$ and a smooth $4$-pointed rational
curve $(C, q, p_1, p_2, p_3)$ carrying the line bundle $\mathcal{O}_C$;
\item $\Delta_{1, \{ 1 \}}$, whose general member is obtained by
joining with an ordinary node a smooth $2$-pointed elliptic curve
$(E, p, p_1)$ carrying a square root of $\mathcal{O}_E(p + p_1)$
and a smooth $3$-pointed rational curve $(C, q, p_2, p_3)$
carrying the line bundle $\mathcal{O}_C$;
\item $\Delta_{1, \{ 2 \}}$, whose general member is obtained by
joining with an ordinary node a smooth $2$-pointed elliptic curve
$(E, p, p_2)$ carrying a square root of $\mathcal{O}_E(p + p_2)$
and a smooth $3$-pointed rational curve $(C, q, p_1, p_3)$
carrying the line bundle $\mathcal{O}_C$;
\item $A_{1, \{ 3 \}}$, whose general member is obtained by
joining with an exceptional component a smooth $2$-pointed elliptic curve
$(E, p, p_3)$ carrying an even root of $\mathcal{O}_E$ and a smooth
$3$-pointed rational curve $(C, q, p_1, p_2)$ carrying the line bundle $\mathcal{O}_C$;
\item $B_{1, \{ 3 \}}$, whose general member is obtained by
joining with an exceptional component a smooth $2$-pointed elliptic curve $(E, p, p_3)$
carrying the line bundle $\mathcal{O}_E$ and a smooth $3$-pointed rational
curve $(C, q, p_1, p_2)$ carrying the line bundle $\mathcal{O}_C$.
\end{itemize}

\begin{Lemma}\label{h^2bis}
The vector space $H^2(\overline{S}_{1,2}^{\hspace{0.05cm}(1,1)})$ is
generated by boundary classes.
\end{Lemma}

\proof We claim that the four boundary classes $\alpha_\irr$,
$\beta_\irr$, $\alpha_{1, \emptyset}$, and $\beta_{1,\emptyset}$ are
linearly independent. Indeed, suppose that there is a relation:
\begin{equation}\label{relation}
a_0 \alpha_\irr + b_0 \beta_\irr + a_1 \alpha_{1, \emptyset} +
b_1 \beta_{1, \emptyset} = 0.
\end{equation}
By restricting (\ref{relation}) to $B_\irr$, $A_\irr$,
$B_{1,\emptyset}$ and $A_{1,\emptyset}$, we obtain $a_1=0$, $b_1=0$,
$a_0=0$ and $b_0=0$, respectively. Since we already know that the
first Betti number vanishes, it will be sufficient to prove that
$$
\chi(\overline{S}_{1,2}^{\hspace{0.05cm}(1,1)})=6.
$$
We first compute $\chi(S_{1,2}^{\hspace{0.05cm}(1,1)})$. A point of
$S_{1,2}^{\hspace{0.05cm}(1,1)}$ corresponds to a smooth $2$-pointed
elliptic curve $(E; p_1,p_2)$ together with a square root of the
line bundle $\mathcal{O}_E(p_1+p_2)$, that is, a ramification point
of the $2$-sheeted covering of $\bP^1$ defined by the linear series
$\vert p_1 + p_2 \vert$. Hence the natural projection
$S_{1,2}^{\hspace{0.05cm}(1,1)} \to M_{1,2}$ is generically
four-to-one, but there are a few special fibers with less than four
points. Indeed, consider the two-sheeted covering of $\bP^1$ defined
by the linear series $\vert 2 p_1 \vert$ and ramified over $\infty$,
$0$, $1$, and $\lambda$, with $p_1$ lying above $\infty$. If $p_2$
lies above $0$, then the corresponding involution exchanges cyclically
the square roots of $\mathcal{O}_E(p_1+p_2)$. If moreover $\lambda =
- 1$, then the projectivity of $\bP^1$ defined by $z \mapsto - z$
induces another automorphism of $(E; p_1, p_2)$ and in this case all
square roots of $\mathcal{O}_E(p_1+p_2)$ are identified. Finally, if
$\lambda = - \omega$ (with $\omega^3 = 1$) and $p_2$ is one point
lying above $\frac{\omega}{\omega - 1}$ then the projectivity of
$\bP^1$ defined by $z \mapsto \frac{z+\omega}{\omega}$ induces
automorphisms of $(E, p_1, p_2)$ exchanging ciclically three square
roots of $\mathcal{O}_E(p_1+p_2)$. As in \cite{ArbCor:98}, p.~124,
we denote by $X$ the locus of all curves $(E; p_1,p_2)$ such that
$p_2$ is a $2$-torsion point with respect to the group law with
origin in $p_1$. Therefore, since $X$ is isomorphic to the quotient
$M'_{0,4}$ of $M_{0,4}$ modulo the operation of interchanging the
labelling of two of the marked points, we have:
\begin{eqnarray*}
\chi(S_{1,2}^{\hspace{0.05cm}(1,1)}) &=& 4 \chi(M_{1,2} \setminus X
\cup \{ \textrm{point} \})
+ 2 \chi(X \setminus \{ \textrm{point} \}) + \chi(\textrm{point}) + \\
& & + 2 \chi(\textrm{point}) = 4 \chi(M_{1,2}) - 2 \chi(M'_{0,4}) -
3 = 1
\end{eqnarray*}
(recall that $\chi(M'_{0,4})=0$ and $\chi(M_{1,2}) = 1$ by
\cite{ArbCor:98}, (5.3) and (5.4)).

Next, from the stratification of $\overline{M}_{1,2}$ by graph type
(see \cite{ArbCor:98}, Figure~1), it follows that
\begin{eqnarray*}
\chi(\overline{S}_{1,2}^{\hspace{0.05cm}(1,1)}) =
\chi(S_{1,2}^{\hspace{0.05cm}(1,1)}) + 2 \chi(M'_{0,4}) +
\chi(S_{1,1}^{\hspace{0.05cm}(0),+}) +
\chi(S_{1,1}^{\hspace{0.05cm}(0),-}) + 3 + 1.
\end{eqnarray*}
Since $\chi(S_{1,1}^{\hspace{0.05cm}(0),+}) = 0$ by \cite{BF:09}, (5),
and $\chi(S_{1,1}^{\hspace{0.05cm}(0),-}) = \chi(M_{1,1}) = 1$,
we obtain $\chi(\overline{S}_{1,2}^{\hspace{0.05cm}(1,1)}) = 6$,
as claimed.

\qed

The following result is a partial analogue of \cite{BF:09}, Lemma~3.

\begin{Lemma}\label{kernelbis}
Let $x$ and $y$ be distinct and not belonging to $\{ 1,2,3 \}$. Let
${\mathfrak S}_2$ be the symmetric group permuting $x$ and $y$.
Define
$$
\xi: \overline{M}_{0, \{ 1, 2, 3, x, y \}}/{\mathfrak S}_2
\longrightarrow B_\irr \hookrightarrow
\overline{S}_{1,3}^{\hspace{0.05cm}(1,1,0)}
$$
by joining the points labelled $x$ and $y$ with an exceptional component.
Then the kernel of
$$
\xi^*: H^2(\overline{S}_{1,3}^{\hspace{0.05cm}(1,1,0)})
\longrightarrow  H^2(\overline{M}_{0, \{ 1,2,3,x,y \}}/{\mathfrak
S}_2)
$$
is four-dimensional and generated by $\alpha_\irr$, $\beta_\irr$,
$\beta_{1, \emptyset}$, and $\beta_{1, \{ 3 \}}$.
\end{Lemma}

\proof It is clear that $\xi^*(\alpha_\irr) = \xi^*(\beta_{1, \emptyset})
= \xi^*(\beta_{1, \{ 3 \}})= 0$ because all the corresponding boundary
divisors are disjoint from $B_\irr$. Next, from $\xi^*(\alpha_\irr)=0$ and
\cite{ArbCor:98}, Lemma 3.16, it follows that also $\xi^*(\beta_\irr) = 0$.
Conversely, if $\xi^*(\alpha)=0$ then we claim that for a suitable choice
of rational coefficients $x,y,z,w$ the class $\gamma = \alpha - x \alpha_\irr
- y \beta_\irr - z \beta_{1, \emptyset} - w \beta_{1, \{ 3 \}}$ vanishes on
$\overline{S}_{1,3}^{\hspace{0.05cm}(1,1,0)}$, hence the class $\alpha$
is a linear combination of $\alpha_\irr$, $\beta_\irr$, $\beta_{1, \emptyset}$,
and $\beta_{1, \{ 3 \}}$. In order to show that $\gamma = 0$ we first check
that its restriction to all boundary components vanishes and then we observe
that the restriction map is injective in our case. More precisely, we define
$$
\zeta: \overline{M}_{0, \{ 1, 2, 3, x, y \}}/{\mathfrak S}_2
\longrightarrow A_\irr \hookrightarrow
\overline{S}_{1,3}^{\hspace{0.05cm}(1,1,0)}
$$
by joining the points labelled $x$ and $y$ with an ordinary node.

\noindent Since $h^2(\overline{M}_{0, \{ 1, 2, 3, x, y
\}}/{\mathfrak S}_2)=4$, for any class $\alpha \in
H^2(\overline{S}_{1,3}^{\hspace{0.05cm}(1,1,0)})$ we have
$$
\zeta^*(\alpha)= h \zeta^*(\beta_{1, \emptyset}) + k
\zeta^*(\beta_{1, \{ 3 \}}) + s \zeta^*(\delta_{1, \{ 1 \}}) + t
\zeta^*(\delta_{1, \{ 2 \}})
$$
for some $h,k,s,t \in {\mathbb Q}$.

Suppose now $\xi^*(\alpha)=0$. Then the class
$\beta := \alpha - h \beta_{1, \emptyset} - k \beta_{1, \{ 3 \}}$
satisfies $\zeta^*(\beta)= s \zeta^*(\delta_{1, \{ 1 \}}) + t
\zeta^*(\delta_{1, \{ 2 \}})$ and $\xi^*(\beta)=0$. If
\begin{eqnarray*}
\rho: H^2(\overline{S}_{1,3}^{\hspace{0.05cm}(1,1,0)}) &\longrightarrow&
H^2(A_\irr) \oplus H^2(B_\irr) \oplus H^2(A_{1, \emptyset})
\oplus H^2(B_{1, \emptyset}) \oplus \\
& & H^2(\Delta_{1, \{ 1 \}}) \oplus H^2(\Delta_{1, \{ 2 \}}) \oplus
H^2(A_{1, \{ 3 \}}) \oplus H^2(B_{1, \{ 3 \}})
\end{eqnarray*}
is the restriction to the boundary components, we have:
\begin{eqnarray*}
\rho (\beta) &=&
(\zeta^*(\beta), 0, (a \delta_\irr, 0), \beta.B_{1, \emptyset},
c_1 \alpha_\irr + c_2 \beta_\irr + c_3 \alpha_{1, \emptyset} + c_4 \beta_{1, \emptyset}, \\
& & d_1 \alpha_\irr + d_2 \beta_\irr + d_3 \alpha_{1, \emptyset} + d_4 \beta_{1, \emptyset},
e \delta_\irr, \beta.B_{1, \{ 3 \}}),
\end{eqnarray*}
where cohomology classes are expressed in standard bases for the second cohomology groups of the moduli
spaces dominating the various boundary components. In particular, we have $a \delta_\irr \in
H^2(\overline{S}_{1,1}^{\hspace{0.05cm}(0),+})$, $c_1 \alpha_\irr + c_2 \beta_\irr + c_3 \alpha_{1, \emptyset}
+ c_4 \beta_{1, \emptyset} \in H^2(\overline{S}_{1,2}^{\hspace{0.05cm}(1,1)})$,
$d_1 \alpha_\irr + d_2 \beta_\irr + d_3 \alpha_{1, \emptyset} + d_4 \beta_{1, \emptyset}
\in H^2(\overline{S}_{1,2}^{\hspace{0.05cm}(1,1)})$, $e \delta_\irr \in H^2(\overline{S}_{1,2}^{\hspace{0.05cm}(0,0),+})$.
The vanishing of several coefficients is due to the fact that all the above classes restrict to zero on
$B_\irr$ since $\xi^*(\alpha)=0$. For instance, from Lemma~\ref{h^2bis} it follows that
$H^2(\overline{S}_{1,2}^{\hspace{0.05cm}(1,1)})$ is generated by $\alpha_\irr$, $\beta_\irr$,
$\alpha_{1, \emptyset}$, and $\beta_{1, \emptyset}$. Under the corresponding morphism
$$
\xi': \overline{M}_{0, \{ 1,2,x,y \}}/ {\mathfrak S}_2
\longrightarrow B_\irr \hookrightarrow
\overline{S}_{1,2}^{\hspace{0.05cm}(1,1)}
$$
the class $\alpha_{1, \emptyset}$ pulls back to $\delta_{ 0, \{ x,y \}}$, which is not zero, so the kernel of $\xi'^*$
is generated by $\alpha_\irr$, $\beta_\irr$, and $\beta_{1, \emptyset}$, and we have $c_3=d_3=0$.

Moreover, the various restrictions have to coincide also on all the other overlaps between boundary components.
In particular, a careful case-by-case direct inspection shows the following implications:
\begin{itemize}
\item $A_\irr \cap \Delta_{1, \{ 1 \}} \ne \emptyset
\Rightarrow s=c_4=c_3=0$;
\item $A_\irr \cap \Delta_{1, \{ 2 \}} \ne \emptyset
\Rightarrow t=d_4=d_3=0$;
\item $B_{1, \emptyset} \cap A_\irr \ne \emptyset
\Rightarrow \beta.B_{1, \emptyset} = (b \delta_\irr, 0)$;
\item $B_{1, \{ 3 \}} \cap A_\irr \ne \emptyset
\Rightarrow \beta.B_{1, \{ 3 \}} = f \delta_\irr$;
\item $A_{1, \{ 3 \}} \cap A_{1, \emptyset} \ne \emptyset
\Rightarrow e = a$;
\item $B_{1, \{ 3 \}} \cap B_{1, \emptyset} \ne \emptyset
\Rightarrow f = b$;
\item $\Delta_{1, \{ 1 \}} \cap A_{1, \emptyset} \ne \emptyset
\Rightarrow c_1+c_2 = a$;
\item $\Delta_{1, \{ 1 \}} \cap B_{1, \emptyset} \ne \emptyset
\Rightarrow c_1 = b$;
\item $\Delta_{1, \{ 2 \}} \cap A_{1, \emptyset} \ne \emptyset
\Rightarrow d_1+d_2 = a$;
\item $\Delta_{1, \{ 2 \}} \cap B_{1, \emptyset} \ne \emptyset
\Rightarrow d_1 = b$.
\end{itemize}
As a consequence, if
$$
\gamma := \beta - b \alpha_\irr - (a-b) \beta_\irr = \alpha - b \alpha_\irr - (a-b) \beta_\irr
- h \beta_{1, \emptyset} - k \beta_{1, \{ 3 \}}
$$
then
$$
\rho(\gamma) = (0, 0, (0,0), (0, 0), 0, 0, 0, 0).
$$
On the other hand, Proposition~\ref{harer} with $k=4, g=1, n=3$, and the inductive argument
following (\ref{exact}) imply that $\rho$ is injective, hence $\gamma=0$ and our claim is proved.

\qed

\begin{Lemma}\label{boundary}
The vector space $H^2(\overline{S}_{1,3}^{\hspace{0.05cm}(1,1,0)})$
is generated by boundary classes.
\end{Lemma}

\proof Let $V$ the subspace of
$H^2(\overline{S}_{1,3}^{\hspace{0.05cm}(1,1,0)})$
generated by the pull-backs $\delta_{1, \emptyset}$,
$\delta_{1, \{ 1 \}}$, $\delta_{1, \{ 2 \}}$, and
$\delta_{1, \{ 3 \}}$ of the corresponding boundary
divisors on $\overline{M}_{1,3}$. In view of Lemma~\ref{kernelbis},
it will be sufficient to show that the morphism $\xi^*$
vanishes modulo $V$. In order to do so,
we adapt the argument in \cite{ArbCor:98}, pp.~114--118.
If $\alpha$ is any class in
$H^2(\overline{S}_{1,3}^{\hspace{0.05cm}(1,1,0)})$, then
$$
\xi^*(\alpha) = a_{ \{x,y \}} \delta_{ 0, \{ x,y \}} +
a_{ \{1,3 \}} \delta_{ 0, \{ 1,3 \}} + a_{ \{2,3 \}} \delta_{ 0, \{ 2,3 \}} +
a_{ \{ 3 \}} ( \delta_{ 0, \{ x,3 \}} + \delta_{ 0, \{ y,3 \}} ).
$$
The idea is simply to modify $\alpha$ with elements of $V$ in such a way
that $\xi^*(\alpha)=0$.
The first move consists in adding to $\alpha$ a suitable multiple of
$\delta_{1, \emptyset}$ so as to make $a_{ \{x,y \}} = 0$.
Next, by Lemma~\ref{h^2bis} there exists a suitable multiple of
$\delta_{1, \{ 3 \}}$ which added to $\alpha$ makes $a_{ \{ 3 \}} = 0$
(for details, see \cite{ArbCor:98}, p.~116). Finally, the third move
consists in adding to $\alpha$ a suitable linear combination
of $\delta_{1, \{ 1 \}}$ and $\delta_{1, \{ 2 \}}$ so as to make
$a_{ \{1,3 \}} = a_{ \{2,3 \}} = 0$. Hence we obtain $\xi^*(\alpha)=0$,
as desired.

\qed

Finally, we conclude the proof of (\ref{eq_six}).

\begin{Lemma}
We have $H^3(\overline{S}_{1,3}^{\hspace{0.05cm}(1,1,0)})=0$.
\end{Lemma}

\proof By Lemma~\ref{boundary},
$H^2(\overline{S}_{1,3}^{\hspace{0.05cm}(1,1,0)})$
is generated by the eight boundary classes
$\alpha_\irr$, $\beta_\irr$, $\alpha_{1, \emptyset}$,
$\beta_{1, \emptyset}$,
$\delta_{1, \{ 1 \} }$, $\delta_{1, \{ 2 \} }$,
$\alpha_{1, \{ 3 \} }$, and $\beta_{1, \{ 3 \} }$.
Hence $h^2(\overline{S}_{1,3}^{\hspace{0.05cm}(1,1,0)}) \le 8$;
next, we claim that
$$
\chi(\overline{S}_{1,3}^{\hspace{0.05cm}(1,1,0)}) = 18,
$$
hence $h^3(\overline{S}_{1,3}^{\hspace{0.05cm}(1,1,0)})=0$.

First of all, we compute $\chi(S_{1,3}^{\hspace{0.05cm}(1,1,0)})$.
The natural projection $S_{1,3}^{\hspace{0.05cm}(1,1,0)} \to
M_{1,3}$ is generically four-to-one, but there are a few special
fibers with less than four points. Indeed, denote by $Y$ the locus
of all curves $(E; p_1,p_2,p_3)$ such that both $p_2$ and $p_3$ are
$2$-torsion points with respect to the group law with origin in
$p_1$. Then it is clear that the fiber of $p$ over $Y$ consists of two
points; moreover, if $(E; p_1,p_2,p_3)$ is the $2$-sheeted covering
of $\bP^1$ ramified over $\infty$, $0$, $1$, and $- \omega$ (with
$\omega^3 = 1$), with $p_1$ lying above $\infty$ and $p_2$, $p_3$
lying above $\frac{\omega}{\omega - 1}$, then the projectivity of
$\bP^1$ defined by $z \mapsto \frac{z+\omega}{\omega}$ induces
automorphisms of $(E; p_1, p_2, p_3)$ exchanging cyclically three
square roots of $\mathcal{O}_E(p_2+p_3)$. Next, we claim that
$\chi(Y)=0$. Indeed, it is clear that $M_{0,4} \setminus \{
\textrm{point} \}$ is a $2$-sheeted covering of $Y \setminus \{
\textrm{point} \}$. Hence
$$
\chi(Y)= \frac{\chi(M_{0,4} \setminus  \{ \textrm{point} \})}{2} + 1
= 0,
$$
as claimed. It follows that
$$
\chi(S_{1,3}^{\hspace{0.05cm}(1,1,0)}) = 4 \chi(M_{1,3} \setminus Y
\cup \{ \textrm{point} \}) + 2 \chi(Y) + 2 \chi(\textrm{point}) = -
2
$$
(recall that $\chi(M_{1,3})=0$ by \cite{ArbCor:98}, (5.4)).

Next, from the stratification of $\overline{M}_{1,3}$ by graph type
(see \cite{ArbCor:98}, Figure~2) we obtain
\begin{eqnarray*}
\chi(\overline{S}_{1,3}^{\hspace{0.05cm}(1,1,0)}) &=&
\chi(S_{1,3}^{\hspace{0.05cm}(1,1,0)}) + 2 \chi(M'_{0,5}) +
\chi(S_{1,1}^{\hspace{0.05cm}(0),+})\chi(M_{0,4}) + \\
& & + \chi(S_{1,1}^{\hspace{0.05cm}(0),-})\chi(M_{0,4}) + 2
\chi(S_{1,2}^{\hspace{0.05cm}(1,1)}) +
\chi(S_{1,2}^{\hspace{0.05cm}(0,0),+}) + \\
& & + \chi(S_{1,2}^{\hspace{0.05cm}(0,0),-}) + 3 \chi(M_{0,4}) + 12
\chi(M'_{0,4}) +
3 \chi(S_{1,1}^{\hspace{0.05cm}(0),+}) + \\
& & + 3 \chi(S_{1,1}^{\hspace{0.05cm}(0),-}) +
9 + 5 + 2.
\end{eqnarray*}
Here $M'_{0,4}$ and $M'_{0,5}$ denote the quotient of $M_{0,4}$ and $M_{0,5}$,
respectively, modulo the operation of interchanging the labelling of two of
the marked points; by \cite{ArbCor:98}, (5.3), we have $\chi(M'_{0,4})=0$
and $\chi(M'_{0,5})=1$. Moreover, we have $\chi(S_{1,1}^{\hspace{0.05cm}(0),+})=0$,
$\chi(S_{1,1}^{\hspace{0.05cm}(0),-})=1$, $\chi(S_{1,2}^{\hspace{0.05cm}(1,1)})=1$
(see above, proof of Lemma \ref{h^2bis}), $\chi(S_{1,2}^{\hspace{0.05cm}(0,0),+})=0$
by \cite{BF:09}, (6), $\chi(S_{1,2}^{\hspace{0.05cm}(0,0),-})=\chi(M_{1,2})=1$
by \cite{ArbCor:98}, (5.4), $\chi(M_{0,4})=-1$, hence
$\chi(\overline{S}_{1,3}^{\hspace{0.05cm}(1,1,0)}) = 18$, as desired.

\qed

\section{The Picard group and the second cohomology group}

Let $\Gamma$ be a finite index subgroup of the mapping class group
$\Gamma_{g,n}$ that contains the Torelli group. For instance,
$\Gamma$ can be any of the levels in Example A or Example B in
Section~\ref{sec2}.

We recall that the Hodge class is defined to be the first Chern
class $\lambda:=c_1(\pi_*\omega_{\pi})$, where $\omega_{\pi}$ is the
relative dualizing sheaf of the universal curve $\pi$. For any
level $\Gamma$ denote by $p_{\Gamma}: {\overline M}^{\Gamma}
\rightarrow {\overline M}_{g,n}$ the map induced by the definition
of ${\overline M}^{\Gamma}$. Define on ${\overline M}^{\Gamma}$ the
Hodge class to be the pull-back $p_{\Gamma}^*(\lambda)$, which we
still denote by $\lambda$ by abuse of notation.

The universal cotangent classes of ${\overline M}_{g,n}$ are the
first Chern classes of the line bundles ${\mathcal L}_i$ for $i=1,
\ldots, n$, where the fiber of ${\mathcal L}_i$ over a pointed
stable curve $[C;x_1, \ldots, x_n]$ is the cotangent space
$T^*_{x_i}(C)$ - we recall that the marked points of any stable
curve are smooth. Naturally, the pull-backs under $p_{\Gamma}^*$ of
the classes $\psi_i$'s define universal cotangent classes on any
level structure $M^{\Gamma}$.

Finally, the irreducible components of the Deligne-Mumford boundary
of a level structure yield degree two cohomology classes via
(rational) Poincar\'{e} duality. We refer to them as {\em boundary
classes}.

This said, we can prove the following theorem.

\begin{Theorem}\label{pic}
Let $\Gamma$ be a finite-index subgroup of the mapping class group
that contains the Torelli group.  Let $\overline{M}^{\Gamma}$ be the
corresponding covering over $\overline{M}_{g,n}$. If $g \ge 5$, then
$\Pic (\overline{M}^\Gamma) \otimes \mathbb{Q}$ is freely generated
by the Hodge class, the set of $\psi$ classes, and the set of
boundary classes.
\end{Theorem}

\proof From \cite{Put}, Theorem 2.1, it follows that $H^2(M^\Gamma)$
is freely generated by the Hodge class and the set of $\psi$
classes. On the other hand, by \cite{Hain}, Theorem 5.4, we have
$H^2(M^\Gamma, \mathbb{Q}) \cong \Pic (M^\Gamma) \otimes \mathbb{Q}$,
hence we need only to check that the boundary classes are linearly
independent. Indeed, by the proof of Theorem 10 in \cite{ArbCor:08}
$H^1(\overline{M}^\Gamma)=0$ implies that $H^1(M^\Gamma)$ is the
kernel of the Gysin map $H^0(\partial \overline{M}^\Gamma) \to
H^2(\overline{M}^\Gamma)$. Since by \cite{Hain}, Proposition 5.2, we
have $H^1(M^\Gamma)=0$, it follows that there are no linear
relations among boundary classes in $H^2(\overline{M}^\Gamma)$,
hence the claim holds. In the special case of spin moduli spaces,
the linear independence of boundary classes can be also checked
directly by intersecting them with suitable test families, see
\cite{Cor}, Proposition (7.2).

\qed

Now we focus on the special case of $\sgnm$. As pointed out in
\cite{BF:06}, the injection (\ref{injective}) is compatible with the
Hodge decomposition. As a consequence, Proposition~\ref{harer} and
Poincar\'e duality imply that $h^{2,0}(\sgnm)=0$ for every $g,n$ if
$h^{2,0}(\overline{S}_{0,n}^{\hspace{0.05cm}(m_1, \ldots, m_n)})=0$,
$n \le 5$, and $h^{2,0}(\overline{S}_{1,n}^{\hspace{0.05cm}(m_1,
\ldots, m_n)})=0$, $n \le 2$. On the other hand,
$H^2(\overline{S}_{0,n}^{\hspace{0.05cm}(m_1, \ldots, m_n)})$ is
algebraic by (\ref{zero}) and \cite{Keel:92}. Moreover,
$H^2(\overline{S}_{1,n}^{\hspace{0.05cm}(0, \ldots, 0)})$ is also
algebraic by \cite{BF:09}, Proposition 2, and the same is true for
$H^2(\overline{S}_{1,2}^{\hspace{0.05cm}(1,1)})$ by
Lemma~\ref{h^2bis}. Hence $h^{2,0}(\sgnm)$ vanishes and from the
exponential sequence
$$
0 \to \Z \to \mathcal{O} \to \mathcal{O}^* \to 0
$$
we deduce that
$$
H^2(\sgnm) \cong H^1(\sgnm, \mathcal{O}^*) = \Pic(\sgnm).
$$

In particular, the following holds.

\begin{Corollary}\label{h2}
If $g \ge 5$ then $H^2(\sgnm)$ is freely generated by the Hodge
class, the set of $\psi$ classes, and the set of boundary classes.
\end{Corollary}

\end{document}